\documentclass[british,review,number]{elsarticle}
\usepackage[T1]{fontenc}
\usepackage[latin9]{inputenc}
\usepackage[active]{srcltx}
\usepackage{color}
\usepackage{amsmath}
\usepackage{amssymb}
\usepackage{graphicx}
\usepackage{esint}

\makeatletter

\providecommand{\tabularnewline}{\\}

\usepackage{epstopdf}
\usepackage{subfigure}\usepackage{upgreek}\usepackage{mathrsfs}\usepackage{multirow}

\journal{Finite Elements in Analysis and Design}

\usepackage{babel}
\usepackage{pdfsync}

\makeatother

\usepackage{babel}
\begin{document}

\title{Towards Automatic Stress Analysis using Scaled Boundary Finite Element
Method with Quadtree Mesh of High-order Elements}

\author{Hou~Man\corref{cor1}}

\ead{h.man@unsw.edu.au}

\cortext[cor1]{Corresponding author. Tel.: +612 93855030}

\author{Chongmin~Song, Sundararajan~Natarajan, Ean~Tat~Ooi, Carolin~Birk\corref{}}

\address{School of Civil and Environmental Engineering}

\address{The University of New South Wales, Sydney, NSW 2052, Australia}
\begin{abstract}
This paper presents a technique for stress and fracture analysis by
using the scaled boundary finite element method (SBFEM) with quadtree
mesh of high-order elements. The cells of the quadtree mesh are modelled
as scaled boundary polygons that can have any number of edges, be
of any high orders and represent the stress singularity around a crack
tip accurately without asymptotic enrichment or other special techniques.
Owing to these features, a simple and automatic meshing algorithm
is devised. No special treatment is required for the hanging nodes
and no displacement incompatibility occurs. Curved boundaries and
cracks are modelled without excessive local refinement. Five numerical
examples are presented to demonstrate the simplicity and applicability
of the proposed technique.\end{abstract}
\maketitle

\section{Introduction}

Finite Element Analysis (FEA) is the most widely used analysing tool
in Computer Aided Engineering (CAE). One key factor to achieve an
accurate FEA is the layout of the finite element mesh, including both
mesh density and element shape \citep{Yerry1983a}. Regions containing
complex boundaries, rapid transitions between geometric features or
singularities require finer discretisation \citep{Cheng1996,Greaves1999}.
This leads to the development of adaptive meshing techniques that
assure the solution accuracy without sacrificing the computational
efficiency \citep{Tabarraei2005a,Lo2010}. The construction of a high
quality mesh, in general, takes the most of the analysis time \citep{Hughes2005}.
The recent rapid development of the isogeometric analysis \citep{Hughes2005,Nguyen-Thanh2011a,Simpson2013},
which suppressed the meshing process, has emphasised the significance
of mesh automation in engineering design and analysis.

\begin{figure}
\noindent \begin{centering}
\includegraphics[scale=0.8]{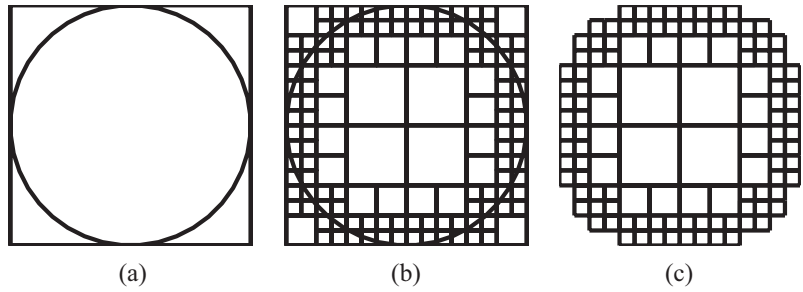}
\par\end{centering}

\caption{\label{fig:qtreerep1}Generation of quadtree mesh on a circular domain.
(a) Cover the domain with a square root cell (b) Subdivide the square
cells (c) Select the cells based on the domain boundary}
\end{figure}

\begin{figure}
\noindent \begin{centering}
\includegraphics[scale=0.5]{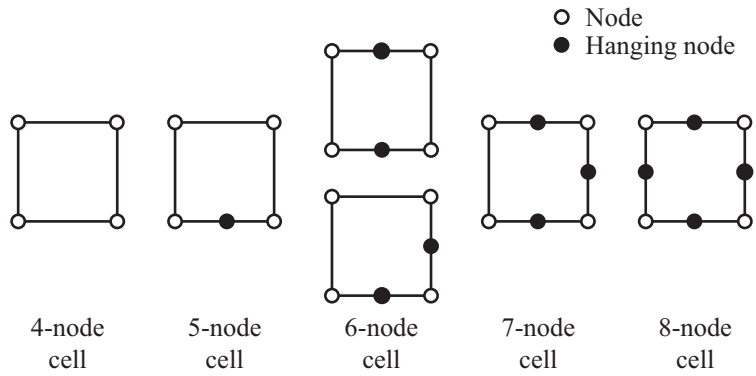}
\par\end{centering}

\caption{\label{fig:qtreecell}6 main types of master quadtree cells with $2:1$
rule enforced}
\end{figure}

Quadtree in FEA is a kind of hierarchical tree-based techniques for
adaptive meshing of a 2D geometry \citep{Greaves1999}. It discretises
the geometry into a number of square cells of different size. The
process is illustrated in Fig.\,\ref{fig:qtreerep1} using a circular
domain. The geometry is first covered with a single square cell, also
known as the root cell of the quadtree (Fig.\,\ref{fig:qtreerep1}a).
As shown in Fig.\,\ref{fig:qtreerep1}b, the root cell is subdivided
into 4 equal-sized square cells and each of the cells is recursively
subdivided to refine the mesh until certain termination criteria are
reached. In this example, a cell is subdivided to better represent
the boundary of the circle and the subdivision stops when the predefined
maximum number of division is reached. The final mesh is obtained
after deleting all the cells outside the domain (Fig.\,\ref{fig:qtreerep1}c).
The cell information is stored in a tree-type data structure, in which
the root cell is at the highest level. It is common practice to limit
the maximum difference of the division levels between two adjacent
cells to one \citep{Yerry1983a,GVS2001}. This is referred to the
$2:1$ rule and the resulting mesh is called a balanced \citep{GVS2001}
or restricted quadtree mesh \citep{Tabarraei2005a}. A balanced quadtree
mesh not only ensures there is no large size difference between adjacent
cells, but also reduces the types of quadtree cells in a mesh to the
6 shown in Fig.\,\ref{fig:qtreecell}. Owing to its simplicity and
large degree of flexibility, the quadtree mesh is also recognised
in large-scale flood/tsunami simulations \citep{Liang2008,Popinet2011}
and image processing \citep{Morvan2007}.


\begin{figure}
\noindent \begin{centering}
\includegraphics[scale=0.5]{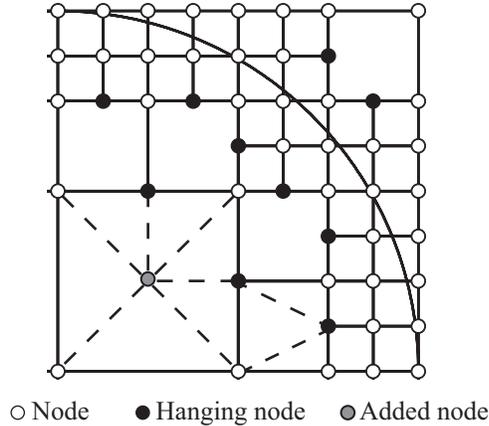}
\par\end{centering}

\caption{\label{fig:qtreerep}Quadtree mesh of the top-right quadrant of a
circular domain. Demonstration of subdivision (dashed lines) is given
in two quadtree cells with hanging nodes on their sides.}
\end{figure}

It is, however, not straightforward to integrate a quadtree mesh in
a FEA. The two major issues are illustrated by Fig.\,\ref{fig:qtreerep},
which shows the quadtree mesh of the top-right quadrant of the circular
domain in Fig.\,\ref{fig:qtreerep1}.
\begin{enumerate}
\item \emph{Hanging nodes} Middle nodes, shown as solid dots in Fig.\,\ref{fig:qtreerep},
exist at the common edges between the adjacent cells with different
division levels. When conventional quadrilateral finite elements are
used, a middle node is connected to the two smaller elements (lower
level) but not to the larger element (higher level). This leads to
incompatible displacement along the edges and the middle nodes are
called the hanging nodes \citep{Greaves1999}. 
\item \emph{Fitting of curved boundary} Quadtree cells are composed of horizontal
and vertical lines only. As shown in Fig.\,\ref{fig:qtreerep}, the
quadtree cells intersected with the curved boundary have to be further
divided into smaller ones to improve the fitting of the boundary.
Generally, the mesh has to be refined in the area surrounding the
boundary. Despite this, the boundary may still not be smooth (Fig.\,\ref{fig:qtreerep1}c)
and may result in unrealistically high stresses. Additional procedure
is required to conform the mesh to the boundary.
\end{enumerate}
There exist a number of different approaches to ensure displacement
compatibility when hanging nodes are present \citep{Ebeida2010,Legrain2011,Tabarraei2008a,Ainsworth2007}.
Three typical approaches among all are briefly discussed here. The
first one is to subdivide the higher level quadtree cells next to
a hanging node into smaller triangular elements \citep{Yerry1983a,Bern1994,Alyavuz2009}
as shown in Fig\,\ref{fig:qtreerep}. Additional nodes may be added
to improve the mesh quality and/or reduce the number of element types.
These techniques lead to a final mesh that only contains conforming
triangular elements. A similar approach was adpoted by \citet{Ebeida2010},
in which the quadtree mesh was subdivided into a conforming mesh dominated
by quadrilateral elements. 

The second approach introduces special conforming shape functions
\citep{Gupta1978} to ensure the displacement compatibility. An early
work by Gupta \citep{Gupta1978} reported the development of a transition
element that had additional node along its side. A conforming set
of shape functions was derived based on the shape functions of the
bilinear quadrilateral elements. Owing to its simplicity and applicability,
Gupta's work was further extended by \citet{Mcdill1987} and \citet{Lo2010}
to hexahedral elements. Fries et al. \citep{Fries2011} investigated
two approaches to handle the hanging nodes within the framework of
the extended finite element method (XFEM). They were different in
whether the enriched degrees-of-freedom (DOFs) were assigned to the
hanging node. A similar work was reported by Legrain et al. \citep{Legrain2011},
in which the selected DOFs are enriched and properly constrained to
ensure the continuity of the field. 

The third approach is to model the quadtree cells as \emph{n-}sided
polygon elements by treating hanging nodes as vertices of the polygon.
This approach generally requires a set of polygonal basis function.
Special techniques are usually required to integrate the resulting
equations over arbitrary polygon domain \citep{Natarajan2009}. This
development was initiated by \citet{wachspress1975rational} who showed
the use of rational basis functions for elements with arbitrary number
of sides. Tabarraei and Sukumar \citep{Tabarraei2005a,Tabarraei2007}
in their work adapted their polygon element \citep{Sukumar2004} to
quadtree mesh. The set of polygonal basis functions was derived using
Laplace interpolant. By using an affine map on the reference polygon,
the conforming shape functions of a quadtree cell with the same number
of vertices (including the hanging nodes) were obtained. They also
reported a fast technique for computing the global stiffness matrix,
making use of the quadtree structure by defining parent elements \citep{Tabarraei2007}.
In this way, the elemental stiffness matrix has to be computed only
15 times (4-node cell not included) for a balanced quadtree mesh (when
$2:1$ rule is enforced). Further development of their work with XFEM
was also reported in \citep{Tabarraei2008a}. 

As mentioned in a recent paper \citep{Sukumar2013}, the development
of high-order polygon elements received relatively less attention.
Mibradt and Pick \citep{Milbradt2008} devised high order basis functions
for polygons based on the natural element coordinates. Those basis
functions, however, are not complete polynomials. Rand et al. \citep{Rand2013}
developed a quadratic serendipity element for arbitrary convex polygons
based on generalised barycentric coordinates. The potential of using
their approach for higher order serendipity elements on convex polygons
was also reported. Based on the same approach, Sukumar \citep{Sukumar2013}
recently developed the quadratic serendipity shape functions that
were applicable for convex and nonconvex polygons and were complete
quadratic polynomials. The shape functions were obtained through solving
an optimisation problem, which was derived from the maximum-entropy
principle. 

Besides dealing with the hanging nodes in a quadtree mesh, fitting
complex boundaries is another challenging part in the mesh generation.
\citet{Yerry1983a} proposed trimming the quadtree cells, intersected
with the boundary, into polygons before a further subdivision process
into triangles or quadrilaterals. Alternatively those cells were first
subdivided and some of the vertices were repositioned based on their
projections onto the boundary \citep{Greaves1999,Alyavuz2009}. In
\citet{Ebeida2010} and \citet{Liang2010}, after the subdivision,
a buffer zone was introduced between the boundary and the internal
quadtree cells. A compatible mesh was then constructed to fill up
this zone. All these techniques require an additional optimisation
step to ensure the final mesh quality. Within the framework of XFEM,
quadtree cells intersected with the boundary were not modified in
pre-processing stage \citep{Fries2011}. However, when constructing
the stiffness matrix, the domain boundary is still required to identify
the portion of the cell within the domain for numerical integration.
In the integration process, that portion of the cell within the domain
is either subdivided into geometric sub-cells \citep{Dreau2010} or
treated as a polygon \citep{Natarajan2010}.

\begin{figure}
\noindent \begin{centering}
\includegraphics[scale=0.7]{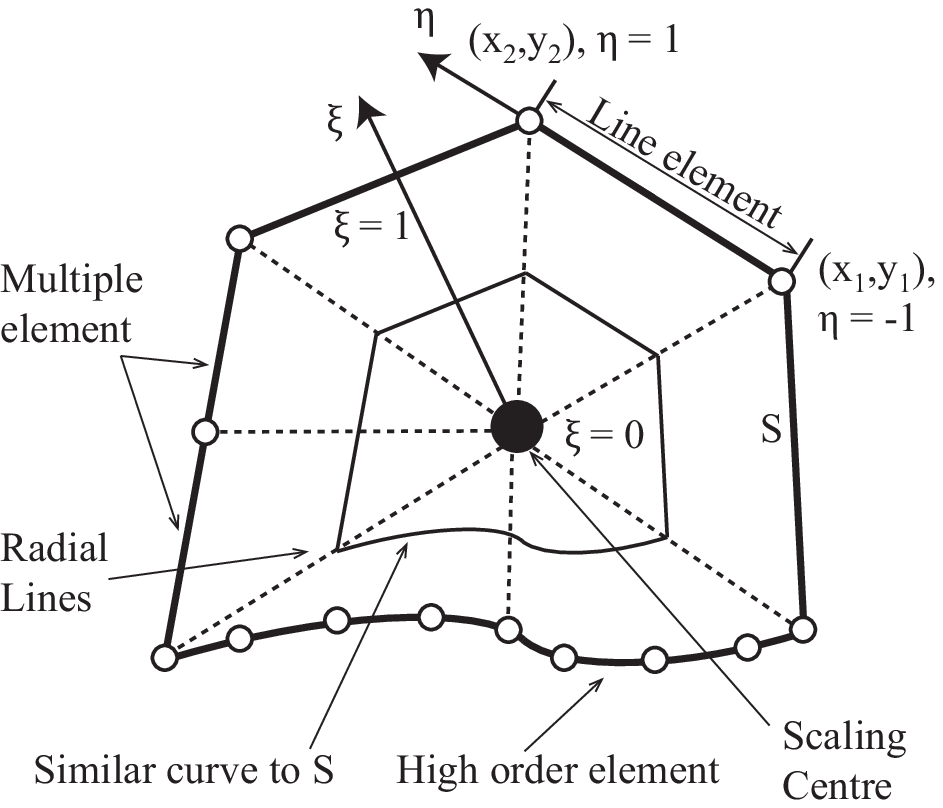}
\par\end{centering}

\caption{\label{fig:sbfempolygon} Scaled boundary representation of a polygon}
\end{figure}

The scaled boundary finite element method (SBFEM) provides an attractive
alternate technique to construct polygon elements (scaled boundary
polygon) \citep{Ooi2012a,Ooi2013} (Fig.\,\ref{fig:sbfempolygon}).
It is a semi-analytical procedure developed by Song and Wolf to solve
boundary value problems \citep{song1997scaled}. The only requirement
for a scaled boundary polygon is that its entire boundary is visible
from the \emph{scaling centre} \citep{song1997scaled}. Only the edges
of the polygon are discretised into line elements. The number of line
elements on an edge can be as many as required. Any type of displacement-based
line elements, including high-order spectral elements\textcolor{black}{{}
\citep{vu2006use}}, can be used. The domain of the scaled boundary
polygon is constructed by scaling from its scaling centre to its boundary,
and the solution within the polygon is expressed semi-analytically
\citep{Ooi2012a,Ooi2013}. A salient feature of the scaled boundary
polygons is that stress singularities occurring at crack and notch
tips, formed by one or several materials, can be accurately modelled
without resorting to asymptotic enrichment and local mesh refinement.
Its high accuracy and flexibility in mesh generation lead to simple
remeshing procedures when modelling crack propagation \citep{yang2006,Ooi2009,Ooi2010a,Ooi2013}. 

\begin{figure}
\noindent \begin{centering}
\includegraphics[scale=0.4]{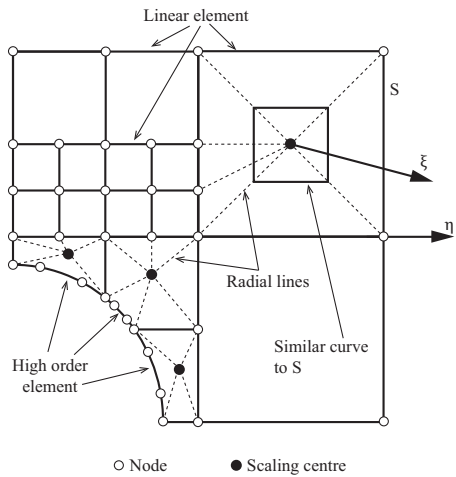}
\par\end{centering}

\caption{\label{fig:sbfequadtreerep}Scaled boundary representation of quadtree
cells}
\end{figure}

This paper presents a technique for the stress and fracture analysis
by integrating the scaled boundary finite element method (SBFEM) with
quadtree mesh of high-order elements. This integrated technique possesses
the following features: 
\begin{enumerate}
\item Hanging nodes are treated without cell subdivision. Each quadtree
cell is modelled as a scaled boundary polygon as shown in Fig.\,\ref{fig:sbfequadtreerep}.
The edges of a quadtree cell can be divided into more than one line
element to ensure displacement compatibility with the adjacent smaller
cells. Hanging-nodes are thus treated the same as other nodes. Owing
to the SBFE formulations \citep{Ooi2012a}, no additional procedure
is required to compute the shape functions for the quadtree cells.
High-order elements can also be used within each quadtree cell directly. 
\item The entire quadtree meshing process is simple and automatic. The boundary
of the problem domain is defined using signed distance functions \citep{Talischi2012}.
Only seed points \citep{Greaves1999} are required to be predefined
to control the mesh density. Owing to the ability of the SBFEM in
constructing polygon elements of, practically, arbitrary shape and
order, the quadtree cells trimmed by curved boundaries are simply
treated as a non-square scaled boundary polygon. High-order elements
can be used to fit curved boundaries closely. The resulting mesh conforms
to the boundary without excessive mesh refinement (see Fig.\,\ref{fig:sbfequadtreerep}). 
\item No local mesh refinement or asymptotic enrichment is required for
a quadtree cell containing a crack tip to accurately model the stress
singularity.
\end{enumerate}
The present paper is organised as follows. The summary of the SBFEM
and its application to quadtree cells are first presented in the next
section. It is followed by the developed algorithm of quadtree mesh
generation in Section\,\ref{sec:Quadtree-mesh-generation}. Five
examples are given in Section\,\ref{sec:Numerical-examples} with
detailed discussion on accuracy and convergence. Finally, conclusions
of the present work are stated in Section\,\ref{sec:Conclusion}.

\section{Scaled boundary finite element method on quadtree cells\label{sec:Scaled-boundary-finite}}

This section summarises the scaled boundary finite element method
for 2D stress and fracture analysis. Only the key equations that are
related to its use with a quadtree mesh are listed. A detailed derivation
of the method based on a virtual work approach is given in \citet{deeks2002virtual}.

\subsection{Element formulation}

The SBFEM can be formulated on quadtree cells by treating each cell
as a polygon with arbitrary number of sides (Fig.~\ref{fig:sbfequadtreerep}).
In each cell, a local coordinate system $(\xi,\eta)$ is defined at
a point called the scaling centre from which the entire boundary is
visible. $\xi$ is the radial coordinate with $\xi=0$ at the scaling
centre and $\xi=1$ at the cell boundary. The edges of each cell are
discretised using one-dimensional finite elements with a local coordinate
$\eta$ having an interval of $-1\leq\eta\leq1$. It is noted that
the hanging nodes appearing in the quadtree structure do not require
any special treatment in the SBFEM formulation. They are simply used
as end nodes of the 1D elements.

The coordinate transformation between the Cartesian $(x,y)$ and the
local $(\xi,\eta)$ coordinate systems are given by the scaled boundary
transformation equations \citep{song1997scaled}:
\begin{align}
\mathbf{x}(\xi,\eta)= & \xi\mathbf{N}(\eta)\mathbf{x}_{\mathrm{b}}\label{eq:coordtrans}
\end{align}
where $\mathbf{x}(\xi,\eta)=[x(\xi,\eta)\; y(\xi,\eta)]^{\mathrm{T}}$
is the Cartesian coordinates of a point in the cell, $\mathbf{N}(\eta)$
is the shape function matrix and $\mathbf{x}_{\mathrm{b}}=[\begin{array}{ccccc}
x_{1} & y_{1} & \ldots & x_{n} & y_{n}\end{array}]^{\mathrm{T}}$ is the vector of nodal coordinates of a cell with $n$ nodes.

The displacement field in each cell $\mathbf{u}(\xi,\eta)$ is interpolated
as
\begin{align}
\mathbf{u}(\xi,\eta)= & \mathbf{N}(\eta)\mathbf{u}(\xi)\label{eq:dispfield}
\end{align}

\noindent where $\mathbf{u}(\xi)$ are radial displacement functions
and are obtained by solving the scaled boundary finite element equation
in displacement \citep{song1997scaled}:
\begin{align}
\mathbf{E}_{0}\xi^{2}\mathbf{u}(\xi)_{,\xi\xi}+(\mathbf{E}_{0}-\mathbf{E}_{1}+\mathbf{E}_{1}^{\mathrm{T}})\xi\mathbf{u}(\xi)_{,\xi}-\mathbf{E}_{2}\mathbf{u}(\xi)= & 0\label{eq:sbfedispeqn}
\end{align}

\noindent with coefficient matrices
\begin{align}
\mathbf{E}_{0}= & \int_{-1}^{+1}\mathrm{\mathbf{B}}_{1}^{\mathrm{T}}(\eta)\mathrm{\mathbf{D\mathrm{\mathbf{B}}_{\mathrm{1}}}(\eta)}|\mathbf{J}(\eta)|d\eta\label{eq:e0}\\
\mathbf{E}_{1}= & \int_{-1}^{+1}\mathrm{\mathbf{B}}_{2}^{\mathrm{T}}(\eta)\mathrm{\mathbf{D\mathrm{\mathbf{B}}_{\mathrm{1}}}(\eta)}|\mathbf{J}(\eta)|d\eta\label{eq:e1}\\
\mathbf{E}_{2}= & \int_{-1}^{+1}\mathrm{\mathbf{B}}_{2}^{\mathrm{T}}(\eta)\mathrm{\mathbf{D\mathrm{\mathbf{B}}_{\mathrm{2}}}(\eta)}|\mathbf{J}(\eta)|d\eta\label{eq:e2}
\end{align}

\noindent where $\mathbf{D}$ is the material constitutive matrix,
$\mathbf{B}_{1}(\eta)$ and $\mathbf{B}_{2}(\eta)$ are the SBFEM
strain-displacement matrices and $|\mathbf{J}(\eta)|$ is the Jacobian
on the boundary required for coordinate transformation.

Eq.\,\eqref{eq:sbfedispeqn} is solved by introducing the variable
$\mathbf{X}(\xi)$ \citep{wolf2003scaled}
\begin{align}
\mathbf{X}(\xi)= & [\begin{array}{cc}
\mathbf{u}(\xi) & \quad\mathbf{q}(\xi)\end{array}]^{\mathrm{T}}\label{eq:Xksi}
\end{align}

\noindent where
\begin{align}
\mathbf{q}(\xi)= & \mathbf{E}_{0}\xi\mathbf{u}(\xi)_{,\xi}+\mathbf{E}_{1}^{\mathrm{T}}\mathbf{u}(\xi)\label{eq:qksi}
\end{align}

\noindent so that Eq.\,\eqref{eq:sbfedispeqn} is transformed into
a first order ordinary differential equation with twice the number
of unknowns:
\begin{align}
\xi\mathbf{X}(\xi)_{,\xi}= & -\mathbf{Z}\mathbf{X}(\xi)\label{eq:firstord}
\end{align}

\noindent with the Hamiltonian matrix $\mathbf{Z}$ \citep{wolf2003scaled}
\begin{align}
\mathbf{Z}= & \left[\begin{array}{cc}
\mathbf{E}_{0}^{-1}\mathbf{E}_{1}^{\mathrm{T}} & -\mathbf{E}_{0}^{-1}\\
\mathbf{E}_{1}\mathbf{E}_{0}^{-1}\mathbf{E}_{1}^{\mathrm{T}}-\mathbf{E}_{2} & -\mathbf{E}_{1}\mathbf{E}_{0}^{-1}
\end{array}\right]\label{eq:hamilton}
\end{align}

An eigenvalue decomposition of the $\mathbf{Z}$ results in
\begin{align}
\mathbf{Z}\left[\begin{array}{cc}
\boldsymbol{\Phi}_{\mathrm{u}}^{\mathrm{(n)}} & \boldsymbol{\Phi}_{\mathrm{u}}^{\mathrm{(p)}}\\
\boldsymbol{\Phi}_{\mathrm{q}}^{\mathrm{(n)}} & \boldsymbol{\Phi}_{\mathrm{q}}^{\mathrm{(p)}}
\end{array}\right]= & \left[\begin{array}{cc}
\boldsymbol{\Phi}_{\mathrm{u}}^{\mathrm{(n)}} & \boldsymbol{\Phi}_{\mathrm{u}}^{\mathrm{(p)}}\\
\boldsymbol{\Phi}_{\mathrm{q}}^{\mathrm{(n)}} & \boldsymbol{\Phi}_{\mathrm{q}}^{\mathrm{(p)}}
\end{array}\right]\left[\begin{array}{cc}
\boldsymbol{\Lambda}^{\mathrm{(n)}} & 0\\
0 & \boldsymbol{\Lambda}^{\mathrm{(p)}}
\end{array}\right]\label{eq:eigendecomp}
\end{align}

\noindent where $\boldsymbol{\Lambda}^{(\mathrm{n})}$ and $\boldsymbol{\Lambda}^{(\mathrm{p})}$
are the eigenvalue matrices with real parts satisfying $\mathrm{Re}(\lambda(\boldsymbol{\Lambda}^{(\mathrm{n})})<0$
and $\mathrm{Re}(\lambda(\boldsymbol{\Lambda}^{(\mathrm{p})})>0$,
respectively. $\boldsymbol{\Phi}_{\mathrm{u}}^{(\mathrm{n})}$ and
$\boldsymbol{\Phi}_{\mathrm{q}}^{(\mathrm{n})}$ are the corresponding
eigenvectors of $\boldsymbol{\Lambda}^{(\mathrm{n})}$ whereas $\boldsymbol{\Phi}_{\mathrm{u}}^{(\mathrm{p})}$
and $\boldsymbol{\Phi}_{\mathrm{q}}^{(\mathrm{p})}$ are the eigenvectors
corresponding to $\boldsymbol{\Lambda}^{(\mathrm{p})}$. For bounded
domains such as those considered in this paper, only the eigenvalues
satisfying $\mathrm{Re}(\lambda(\boldsymbol{\Lambda}^{(\mathrm{n})})<0$
lead to finite displacements at the scaling centre. Using Eq.\,\eqref{eq:eigendecomp}
and Eq.\,\eqref{eq:firstord}, the solutions for $\mathbf{u}(\xi)$
and $\mathbf{q}(\xi)$ are
\begin{align}
\mathbf{u}(\xi)= & \boldsymbol{\Phi}_{\mathrm{u}}^{\mathrm{(n)}}\xi^{-\boldsymbol{\Lambda}^{\mathrm{(n)}}}\mathbf{c}^{\mathrm{(n)}}\label{eq:uksisol}\\
\mathbf{q}(\xi)= & \boldsymbol{\Phi}_{\mathrm{q}}^{\mathrm{(n)}}\xi^{-\boldsymbol{\Lambda}^{\mathrm{(n)}}}\mathbf{c}^{\mathrm{(n)}}\label{eq:qksisol}
\end{align}

The integration constants $\mathbf{c}^{(\mathrm{n})}$ in Eq.\,\eqref{eq:uksisol}
and Eq.\,\eqref{eq:qksisol} are obtained from the nodal displacements
at the cell boundary $\mathbf{u}_{\mathrm{b}}=\mathbf{u}(\xi=1)$
as
\begin{align}
\mathbf{c}^{(\mathrm{n})}= & \left(\boldsymbol{\Phi}_{\mathrm{u}}^{\mathrm{(n)}}\right)^{-1}\mathbf{u}_{\mathrm{b}}\label{eq:intcons}
\end{align}

\noindent The stiffness matrix of each quadtree cell is formulated
as \citep{wolf2003scaled}
\begin{align}
\mathbf{K}= & \boldsymbol{\Phi}_{\mathrm{q}}^{\mathrm{(n)}}\left(\boldsymbol{\Phi}_{\mathrm{u}}^{\mathrm{(n)}}\right)^{-1}\label{eq:stf}
\end{align}

\noindent Substituting Eq.\,\eqref{eq:uksisol} into Eq.\,\eqref{eq:dispfield},
the displacement field in a cell is
\begin{align}
\mathbf{u}(\xi,\eta)= & \mathbf{N}(\eta)\boldsymbol{\Phi}_{\mathrm{u}}^{\mathrm{(n)}}\xi^{-\boldsymbol{\Lambda}^{\mathrm{(n)}}}\mathbf{c}^{\mathrm{(n)}}\label{eq:dispfieldsol}
\end{align}

\noindent Using the Hooke's law and the strain-displacement relationship,
the stress at a point in a cell is \citep{wolf2003scaled}
\begin{align}
\boldsymbol{\sigma}(\xi,\eta)= & \boldsymbol{\Psi}_{\sigma}(\eta)\xi^{-\boldsymbol{\Lambda}^{\mathrm{(n)}}-\mathbf{I}}\mathbf{c}^{\mathrm{(n)}}\label{eq:stresfield}
\end{align}

\noindent where $\boldsymbol{\Psi}_{\sigma}(\eta)=\left[\begin{array}{ccc}
\boldsymbol{\Psi}_{\sigma_{xx}}(\eta) & \boldsymbol{\Psi}_{\sigma_{yy}}(\eta) & \boldsymbol{\Psi}_{\tau_{xy}}(\eta)\end{array}\right]^{\mathrm{T}}$ is the stress mode
\begin{align}
\boldsymbol{\Psi}_{\sigma}(\eta)= & \mathbf{D}\left(-\mathbf{B}_{1}(\eta)\boldsymbol{\Phi}_{\mathrm{u}}^{\mathrm{(n)}}\boldsymbol{\Lambda}^{\mathrm{(n)}}+\mathbf{B}_{2}(\eta)\boldsymbol{\Phi}_{\mathrm{u}}^{\mathrm{(n)}}\right)\label{eq:stresmod}
\end{align}

\subsection{Evaluation of stress intensity factors}

\begin{figure}
\noindent \begin{centering}
\includegraphics[scale=0.4]{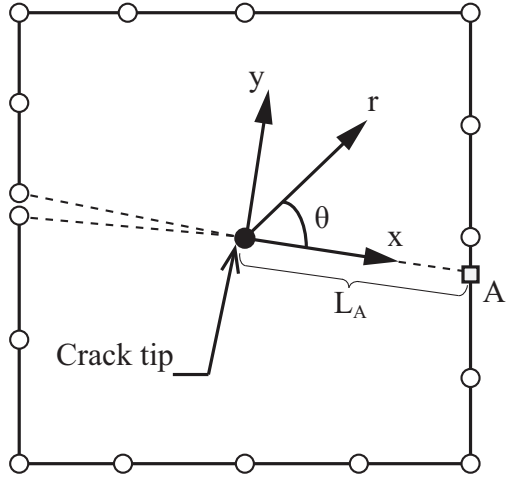}
\par\end{centering}

\caption{\label{fig:crackrep}Modelling of a crack with the scaled boundary
finite element method.}
\end{figure}

Fig.\,\ref{fig:crackrep} shows how a crack is modelled with a quadtree
cell. The crack tip is chosen as the scaling centre. The crack surfaces
are not discretised. The line elements discretising the cell boundary
do not form a closed loop. 

When a crack is modelled by the SBFEM, two eigenvalues, $\lambda_{i}$,
$i=1,\,2$ satisfying $-1<\mathrm{Re}(\lambda_{i})\leq0$ appear in
$\boldsymbol{\Lambda}^{(\mathrm{n})}$. From Eq.\,\ref{eq:stresfield},
it can be discovered that these eigenvalues lead to a stress singularity
as $\xi\rightarrow0$. Using the two modes corresponding to these
two eigenvalues, the singular stresses are expressed as
\begin{align}
\boldsymbol{\sigma}(\xi,\eta)= & \boldsymbol{\Psi}_{\sigma}^{\mathrm{(s)}}(\eta)\xi^{-\boldsymbol{\Lambda}^{\mathrm{(s)}}-\mathbf{I}}\mathbf{c}^{\mathrm{(s)}}\label{eq:singstrefield}
\end{align}

\noindent where
\begin{align}
\boldsymbol{\Lambda}^{(\mathrm{s})}= & \left[\begin{array}{cc}
\lambda_{1} & 0\\
0 & \lambda_{2}
\end{array}\right]\label{eq:singeigenval}
\end{align}

\noindent and $\mathbf{c}^{(\mathrm{s})}$ are the integration constants
corresponding to $\boldsymbol{\Lambda}^{(\mathrm{s})}$. The singular
singular stress modes $\boldsymbol{\Psi}_{\sigma}^{(\mathrm{s})}(\eta)=\left[\begin{array}{ccc}
\boldsymbol{\Psi}_{\sigma_{xx}}^{(\mathrm{s})}(\eta) & \boldsymbol{\Psi}_{\sigma_{yy}}^{(\mathrm{s})}(\eta) & \boldsymbol{\Psi}_{\tau_{xy}}^{(\mathrm{s})}(\eta)\end{array}\right]^{\mathrm{T}}$ is written as
\begin{align}
\boldsymbol{\Psi}_{\sigma}^{\mathrm{(s)}}(\eta)= & \mathbf{D}\left(-\mathbf{B}_{1}(\eta)\boldsymbol{\Phi}_{\mathrm{u}}^{\mathrm{(s)}}\boldsymbol{\Lambda}^{\mathrm{(s)}}+\mathbf{B}_{2}(\eta)\boldsymbol{\Phi}_{\mathrm{u}}^{\mathrm{(s)}}\right)\label{eq:singstremode}
\end{align}

\noindent where $\boldsymbol{\Phi}_{\mathrm{u}}^{(\mathrm{s})}$ are
the modal displacements in $\boldsymbol{\Phi}_{\mathrm{u}}^{(\mathrm{n})}$
corresponding to $\boldsymbol{\Lambda}^{(\mathrm{s})}$.

The stress intensity factors can be computed directly from their definitions.
For a crack that is aligned with the Cartesian coordinate system as
shown in Fig.\,\ref{fig:crackrep}, the stress intensity factors
are defined as
\begin{align}
\left\{ \begin{array}{c}
K_{\mathrm{I}}\\
K_{\mathrm{II}}
\end{array}\right\} = & \lim_{r\rightarrow0}\left\{ \begin{array}{c}
\sqrt{2\pi r}\left.\sigma_{yy}\right|_{\theta=0}\\
\sqrt{2\pi r}\left.\tau_{xy}\right|_{\theta=0}
\end{array}\right\} \label{eq:sifdef}
\end{align}

\noindent Substituting the stress components in Eq.\,\eqref{eq:singstrefield}
into Eq.\,\eqref{eq:sifdef} and using the relation $\xi=r/L_{\mathrm{A}}$
($L_{\mathrm{A}}$ is the distance from the scaling centre to the
boundary along the direction of the crack, see Fig.\,\ref{fig:crackrep})
at $\theta=0$ leads to
\begin{align}
\left\{ \begin{array}{c}
K_{\mathrm{I}}\\
K_{\mathrm{II}}
\end{array}\right\} = & \sqrt{2\pi L_{\mathrm{A}}}\left\{ \begin{array}{c}
\boldsymbol{\Psi}_{\sigma_{yy}}^{\mathrm{(s)}}(\eta(\theta=0))\mathbf{c}^{\mathrm{(s)}}\\
\boldsymbol{\Psi}_{\tau_{xy}}^{\mathrm{(s)}}(\eta(\theta=0))\mathbf{c}^{\mathrm{(s)}}
\end{array}\right\} \label{eq:sifsbfe}
\end{align}

\section{Quadtree mesh generation\label{sec:Quadtree-mesh-generation}}

\begin{figure}
\centering{}\includegraphics[scale=0.8]{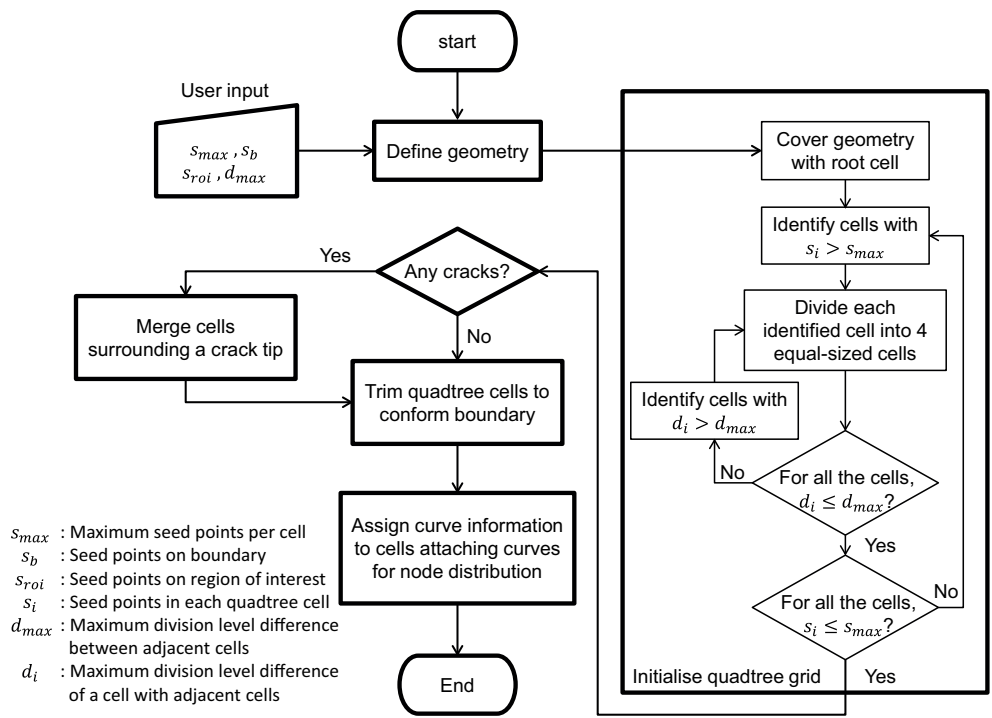}\caption{\label{flowchart}Flow chart of the quadtree mesh generation.}
\end{figure}

\begin{figure}
\centering{}\includegraphics{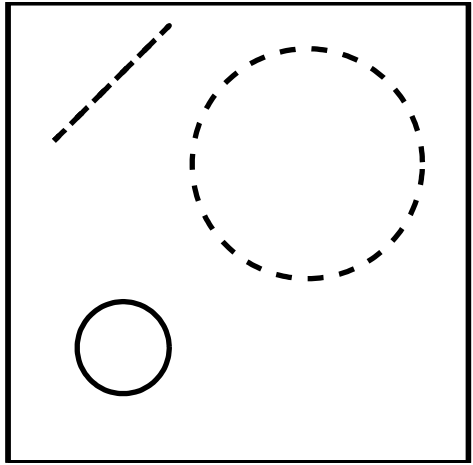}\caption{\label{qtreedes0}Example to illustrate the quadtree mesh generation
process: a square plate with a circular hole. An additional circle
and an inclined line (dashed lines) are included to control local
mesh density.}
\end{figure}

This section presents the developed algorithm for quadtree mesh generation.
Fig.\,\ref{flowchart} shows the flow chart of the overall process.
The entire generation process is automatic with minimal number of
inputs required from the user, which include
\begin{itemize}
\item Maximum allowed number of seed points in a cell $(s_{max})$,
\item Seed points on each boundary $(s_{b})$ and region of interest $(s_{roi})$, 
\item Maximum difference between the division levels of adjacent cells $(d_{max})$,
which is equal to 1 for a balanced quadtree mesh.
\end{itemize}
This section is organised based on Fig\,\ref{flowchart}. It first
presents defining geometry using signed distance function, and assigning
seed points on the boundary and the regions of interest. Detailed
explanations of the meshing steps, which include generating the initial
quadtree grid, trimming the boundary quadtree cells into polygons
and merging cells surrounding a crack tip, are then followed. To facilitate
the description of the meshing steps, Fig.\,\ref{qtreedes0} shows
a square plate with a circular hole and two local refinement features
to be used as an example throughout this section. An efficient computation
of the global stiffness matrix, by taking advantage on the quadtree
mesh, is described at the end of this section. 

\subsection{Define geometry using signed distance function}

\begin{figure}
\noindent \begin{centering}
\includegraphics[scale=0.7]{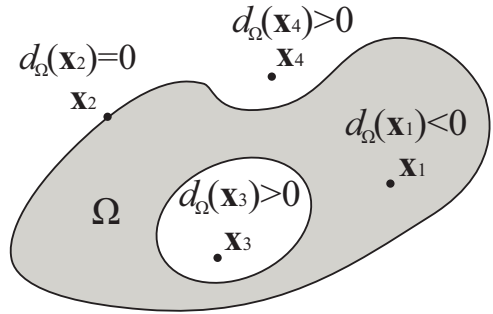}
\par\end{centering}

\caption{\label{dispfunc}Signed distance function of the points inside the
domain ($\mathbf{x}_{1}$), on the boundary ($\mathbf{x}_{2}$) and
outside the domain ($\mathbf{x}_{3}$ and $\mathbf{x}_{4}$)}
\end{figure}

The geometry is defined by using the signed distance function \citep{Persson2004}.
It provides all the essential information of a geometry and can be
operated with simple Boolean operations to build up more complex geometries
\citep{Talischi2012}. The signed distance function of a point $\mathbf{x}\in\mathbb{R}^{2}$
associated with a domain $\Omega$, which is a subset of $\mathbb{R}^{2}$,
is given as

\begin{equation}
d_{\Omega}(\mathbf{x})=s_{\Omega}(\mathbf{x})\min_{\mathbf{y\in\partial\Omega}}\left\Vert \mathbf{x}-\mathbf{y}\right\Vert ,
\end{equation}
where $\partial\Omega$ represents the boundary of the domain and
$\left\Vert \mathbf{x}-\mathbf{y}\right\Vert $ is the \emph{Euclidean
norm} in $\mathbb{R}^{2}$ with $\mathbf{y}\in\partial\Omega$. The
sign function $s_{\Omega}(\mathbf{x})$ is equal $-1$ when $\mathbf{x}$
lies inside the domain and is equal 1 otherwise. This definition of
the signed distance function is visualised in Fig.\,\ref{dispfunc}.
A number of distance functions in MATLAB for simple geometries are
given in \citet{Talischi2012}, including their Boolean operations.

For each boundary and region of interest, a set of pre-defined seed
points \citep{Greaves1999} is introduced to control the quadtree
mesh density. There require four sets of predefined seed points for
the example in Fig.\,\ref{qtreedes0}. Two sets are for the square
and the circular hole representing the actual domain boundary. The
number of seed points directly controls the local density of the quadtree
cells and the quality of fitting the boundary. This is further discussed
in Section\,\ref{sub:Polygon-boundary-cells}. The other two sets
are for the large circle and inclined line controlling local mesh
density only.

\subsection{Initialise quadtree grid}

\begin{figure}
\noindent \begin{centering}
\includegraphics{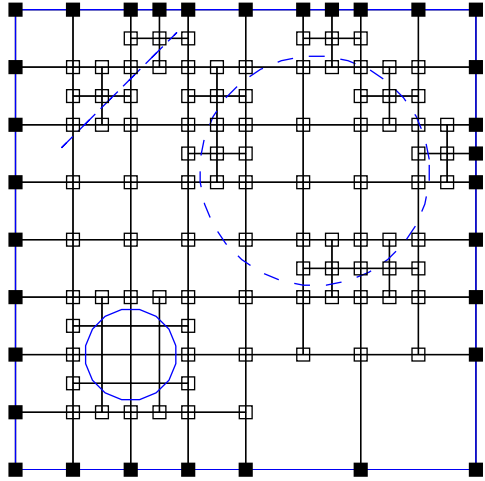}
\par\end{centering}

\caption{\label{qtreedes1}Initial quadtree grid of the example in Fig.~\,\ref{qtreedes0}.
Vertices with solid square markers are on the boundary, with square
box markers are inside the domain, and without any markers are outside
the domain.}
\end{figure}

The meshing process starts with covering the problem domain with a
single square cell (the root cell). The dimension of the root cell
is based on the larger one between the maximum vertical and maximum
horizontal dimension of the geometry. The developed algorithm will
check the number of seed points in the cell. If the number is larger
than the predefined maximum allowed number, the cell will be divided
into 4 equal-sized cells. This generation process is applied recursively
until all the cells have seed points no more than the predefined value.
For each recursive loop, the maximum difference between the division
levels of adjacent cells $(d_{max})$ is enforced. For cells that
have division level difference with the adjacent cells larger than
$d_{max}$, the higher level cell is subdivided into 4 equal-sized
cells. Fig.\,\ref{qtreedes1} shows the initial quadtree grid of
the example in Fig.\,\ref{qtreedes0}. 

\subsection{Trim boundary cells into polygons\label{sub:Polygon-boundary-cells}}

The initial quadtree grid shown in Fig.\,\ref{qtreedes1} does not
conform to the boundary. Those cells that have edges intersected with
the boundary need to be identified and trimmed. By using the signed
distance function, the locations of the vertices (inside the domain,
on the domain boundary or outside the domain as shown in Fig.\,\ref{qtreedes1})
are identified based on the sign and value of the function. For edges
containing two vertices with opposite signs, they are identified as
the edges intersected with the boundary. For each of those edges,
the intersection point with the boundary is computed. 

Some quadtree cells could have vertices very close to the boundary
in comparison with the lengths of their edges. After trimming, poorly
shaped polygon cells with some edges much shorter than the others
could be generated and may adversely affect the mesh quality. To avoid
this situation, the vertices that are within a threshold distance
away from the boundary are identified and then moved to their closest
points on the boundary. In the present work, $1/10$ of the length
of the cell edge (based on the smallest cell attaching to the vertex)
is used as the threshold value. The edges connecting to these vertices
will no longer be cut by the boundary. The trade-off of this process
is the presence of additional non-square cells that lead to additional
computation of the stiffness matrix. This is discussed in Section\,\ref{sub:An-efficient-assembly}.

\begin{figure}
\noindent \begin{centering}
\includegraphics[scale=0.7]{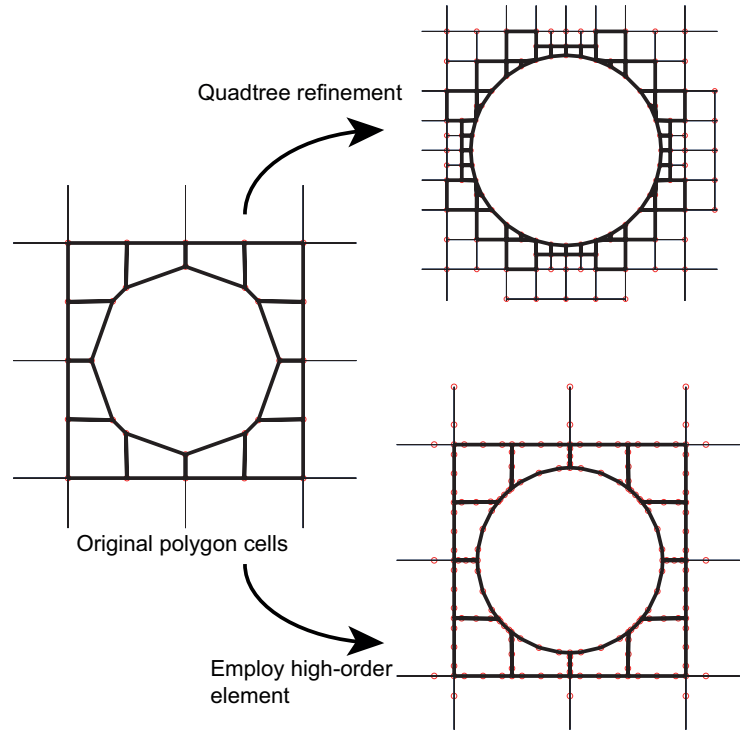}
\par\end{centering}

\caption{\label{qtreedes2}Model curved boundary by quadtree refinement or
using high-order elements. Nodes are represented with small circles
along the cell edges.}
\end{figure}

At the end of the trimming process, the edges of a cell cut by the
boundary are updated with the intersection points and the enclosed
segment of boundary is added to the cell. This will result in polygon
cells. After trimming the quadtree in Fig.\,\ref{qtreedes1}, the
polygon cells around the hole of the example problem is shown in Fig.\,\ref{qtreedes2}.
It is clear from Fig.\,\ref{qtreedes2} that the circular boundary
is not represented accurately if a single linear element is used on
the edge of the cell. 

\begin{figure}
\begin{centering}
\includegraphics{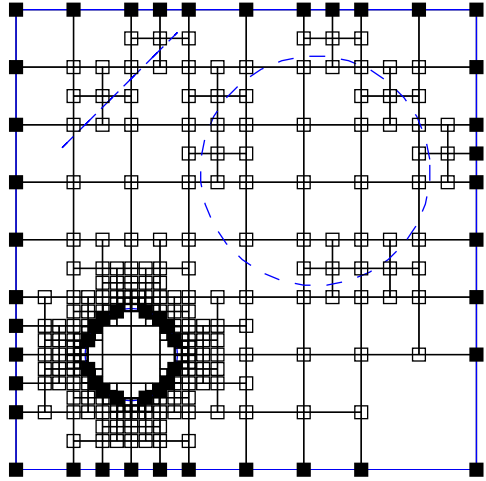}
\par\end{centering}

\caption{\label{qtreedes3}Quadtree mesh after refinement}
\end{figure}

In order to represent the curved boundary more accurately, two alternates
are available in the developed algorithm. The first is to reduce the
element size ($h$-refinement). This is achieved by increasing the
number of seed points on the curved boundary. Fig.\,\ref{qtreedes3}
shows the initial quadtree layout of the example problem after increasing
the seed points around the hole by 4 times. It can be seen by comparing
Fig.\,\ref{qtreedes1} with Fig.~\ref{qtreedes3} that the refinement
is limited to a small region around the hole. The refined quadtree
(Fig\,\ref{qtreedes2}) demonstrates the improvement of capturing
the circular boundary.

The second option to improve the modelling of curved boundaries is
to utilise high-order elements ($p$-refinement). Fig.\,\ref{qtreedes2}
shows the example problem with each line segment on the circular boundary
modelled with a 4th order element. With this approach, curved boundaries
can be captured more accurately using fewer elements. 

Both options to improve the modelling of the boundaries can be applied
simultaneously without conflicts. The numerical accuracy of both approaches
is discussed through numerical examples given in Section\,\ref{sec:Numerical-examples}.

\subsection{Merge cells surrounding a crack tip}

\begin{figure}
\noindent \begin{centering}
\includegraphics[scale=0.7]{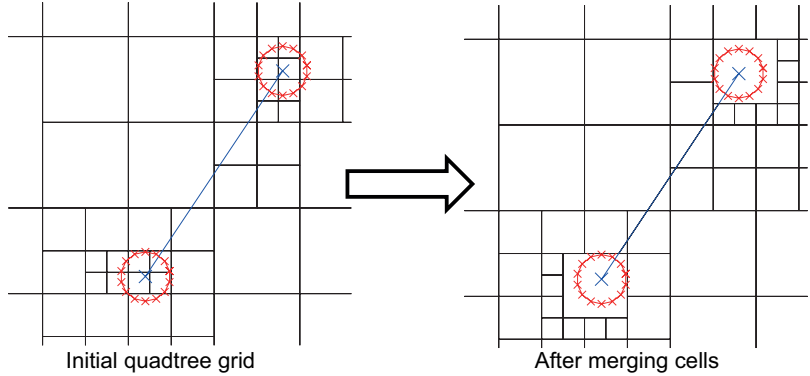}
\par\end{centering}

\caption{\label{qtreedes8}Quadtree mesh for a crack problem before and after
merging cells. The two crack tips are marked with a cross. The two
circles are to control the size of quadtree cells covering the crack
tips.}
\end{figure}

Owing to the capability of the SBFEM for fracture analysis \citep{Song2002},
the domain containing a crack tip is modelled with a single cell.
In the stress solution, the variation along the radial direction,
including the stress singularity, is given analytically and the variation
along the circumference of the cell is represented numerically by
the line elements on boundary. To obtain accurate results, sufficient
nodes have to present on the boundary of the cell to cover the angular
variation of the solution . In the developed algorithm, the size of
a cell containing a crack tip is controlled, as shown in Fig.\,\ref{qtreedes8}
with an inclined crack, by a predefined set of seed points on a circle. 

For problems with cracks, only one additional step is required after
the initial mesh is generated. The cells surrounding the crack tip
are refined to the same division level and then merged into a single
cell as shown in Fig.\,\ref{qtreedes8}. This step avoids having
a crack tip too close to the edges of the cell, which could affect
the mesh quality and the solution accuracy \citep{Ooi2010a}. After
the cells are merged, the intersection point between the edge of the
resulting cell and the crack is computed to define the two crack mouth
points. The other cells on the crack path are split by the crack into
two cells. The splitting process is similar to the trimming of cells
by the boundary, but two vertices are created at every intersection
point between the cell edge and the crack to split the original cells.

\subsection{An efficient construction of the global stiffness matrix\label{sub:An-efficient-assembly}}

The global stiffness matrix is simply the assembly of the stiffness
matrices of each master quadtree and polygon cell. When the $2:1$
rule is enforced to the mesh, only 6 main types of master quadtree
cells are present as given in Fig.\,\ref{fig:qtreecell}. By rotating
the geometry of the master cells orthogonally, the maximum number
of types of these master quadtree cells are 24. For isotropic homogeneous
materials, rotation does not have effect on 4-node or 8-node cells
and only two 2 rotations are required for the first type of 6-node
cell (the top one in Fig.\,\ref{fig:qtreecell}). The maximum number
of master quadtree cells that require stiffness matrix calculation
reduces to 16 (only 15 in \citet{Tabarraei2007} as 4-node cell is
excluded). 

After the mesh generation, the algorithm will check which master cells
out of the 16 appear in the mesh. Their stiffness matrices are then
computed and stored. During the stiffness assembling process, the
stiffness matrix of each regular quadtree cell is directly extracted
from those computed stiffness matrices. For the polygon cells and
those irregular quadtree cells (with their vertices moved to fit the
boundaries), individual stiffness matrix calculation is required.
This approach clearly improves the computational efficiency of constructing
the global stiffness matrix, especially for large scale problems that
contain a significant number of cells. With the use of high-order
elements in the quadtree mesh, this assembling approach becomes even
more economical.

\section{Numerical examples\label{sec:Numerical-examples}}

This section presents five numerical examples to highlight the capability
and the performance of the proposed technique. In the first example,
an infinite plate with a circular hole is modelled and the results
are compared with the analytical solution. The proposed technique
is then used to analyse a square plate with multiple holes to highlight
the automatic meshing capability in handling transition between geometric
features. In the third example, a square plate with a central hole
and multiple cracks is studied to demonstrate the performance of the
proposed technique in handling complicated geometries with singularities.
Thereafter, a square plate with two cracks cross each other is analysed.
It is aimed to emphasise the automation and simplicity of the mesh
generation in the proposed technique. In the first four examples,
the same material properties, with Young's modulus $E=100$ and Poisson's
ratio $v=0.3$, are used. The final example is a cracked nuclear reactor
under internal pressure. It is aimed to show the simplicity of the
present technique in modelling practical non-regular structures.

The computation time reported in this section is based on a desktop
PC with Intel(R) Core(TM) i7 3.40GHz CPU and 16GB of memory. The proposed
technique is implemented in MATLAB and the computation time is extracted
in interactive mode of MATLAB.

\subsection{Infinite plate with a circular hole under uniaxial tension}

\subsubsection{Modelling using exact boundary condition}

\begin{figure}
\begin{centering}
\includegraphics{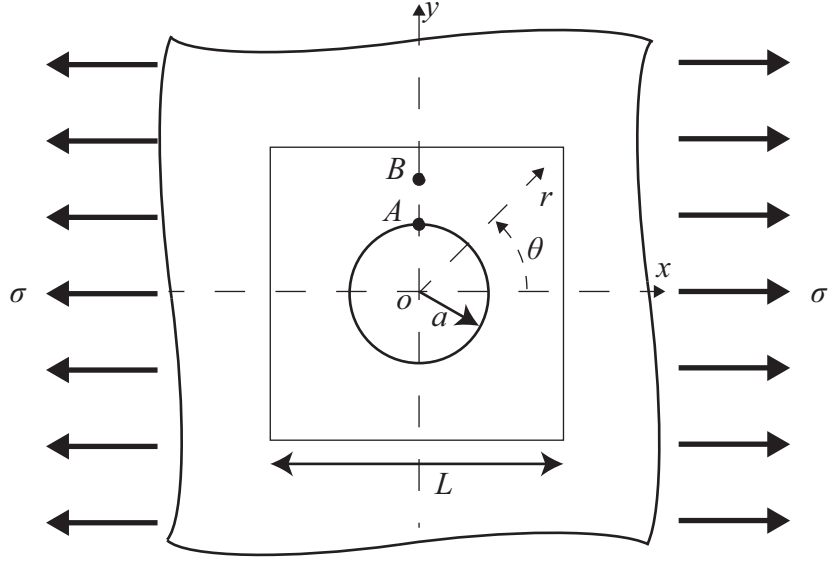}
\par\end{centering}

\caption{Infinite plate with a circular hole under uniaxial tension}
\label{openhole}
\end{figure}

An infinite plate containing a circular hole with radius $a$ at its
centre is considered in this example. The plate is subject to a uniaxial
tensile load as shown in Fig.\,\ref{openhole}. The analytical solution
of the stresses in polar coordinates $(r,\theta)$ is given by \citep{Sukumar2001}:
\begin{align}
\sigma_{11}(r,\theta) & =1-\frac{a^{2}}{r^{2}}\left(\frac{3}{2}\cos2\theta+\cos4\theta\right)+\frac{3a^{4}}{2r^{4}}\cos4\theta\nonumber \\
\sigma_{22}(r,\theta) & =-\frac{a^{2}}{r^{2}}\left(\frac{1}{2}\cos2\theta-\cos4\theta\right)-\frac{3a^{4}}{2r^{4}}\cos4\theta\label{eq:exohstr}\\
\sigma_{12}(r,\theta) & =-\frac{a^{2}}{r^{2}}\left(\frac{1}{2}\sin2\theta+\sin4\theta\right)+\frac{3a^{4}}{2r^{4}}\sin4\theta\nonumber 
\end{align}
The displacement solutions are:
\begin{align}
u_{1}(r,\theta) & =\frac{a}{8\mu}\left[\frac{r}{a}(\kappa+1)\cos\theta+\frac{2a}{r}\left((1+\kappa)\cos\theta+\cos3\theta\right)-\frac{2a^{3}}{r^{3}}\cos3\theta\right]\nonumber \\
u_{2}(r,\theta) & =\frac{a}{8\mu}\left[\frac{r}{a}(\kappa-3)\sin\theta+\frac{2a}{r}\left((1-\kappa)\sin\theta+\sin3\theta\right)-\frac{2a^{3}}{r^{3}}\sin3\theta\right],\label{eq:exohdisp}
\end{align}
where $\mu$ is the shear modulus and $\kappa=\frac{3-v}{1+v}$ is
the Kolosov constant for plane stress condition.

The problem is solved by analysing a finite dimension of the plate
with a dimension of $L\times L$ (see Fig.\,\ref{openhole}). Analytical
traction (Eq.\,\ref{eq:exohstr}) is applied at the four edges of
this finite plate. 

\begin{figure}
\begin{centering}
\subfigure[Global quadtree mesh]{\label{ohm1}\includegraphics{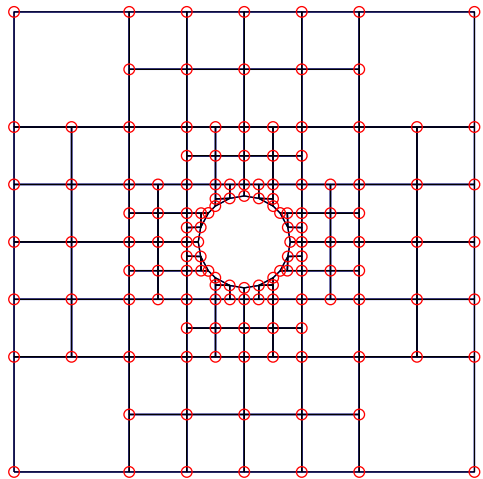}}
\par\end{centering}

\begin{centering}
\subfigure[Cells around the hole]{\label{ohm2}\includegraphics[scale=0.8]{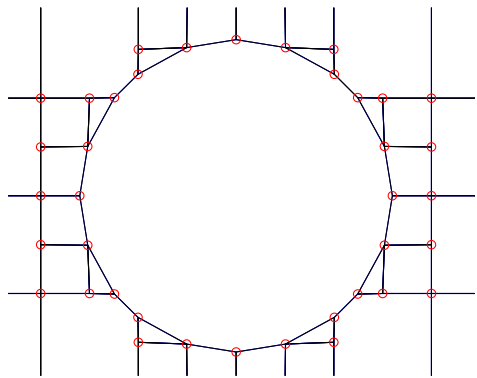}}
\par\end{centering}

\caption{Mesh of a finite square plate with a circular hole ($L/a=10$)}
\label{openholemesh}
\end{figure}

Fig.\,\ref{ohm1} shows the quadtree mesh of the plate for $L/a=10$.
Each edge on a quadtree cell is discretised with 1st order line elements.
The $2:1$ rule is enforced. Based on the proposed technique, the
curved boundary is handled as shown in Fig.\,\ref{ohm2} with polygon
cells. Convergence study is conducted based on the $h-$refinement.
Three different element orders $(p=1,2,4)$ are investigated. 

\begin{figure}
\begin{centering}
\includegraphics[scale=0.8]{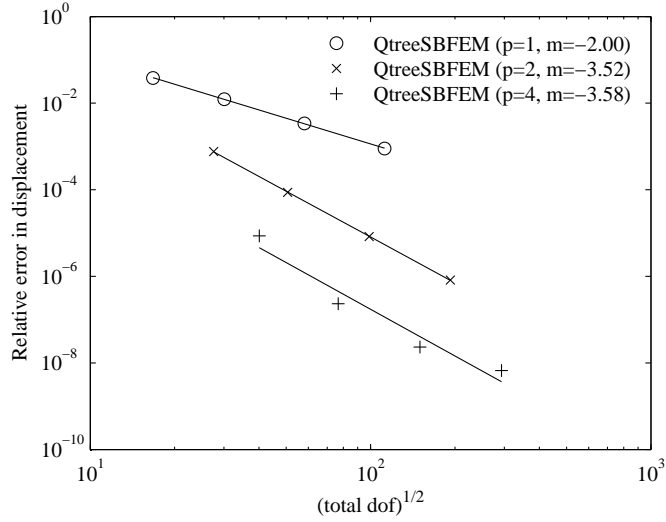}
\par\end{centering}

\caption{Convergence results of the infinite plate with a circular hole, where
$p$ is the element order and $m$ is the slope of the fitted line }
\label{openholecon}
\end{figure}

Fig.\,\ref{openholecon} shows the present results of the relative
error in the displacement norm $\left\Vert {\rm {\rm \mathbf{u}}-{\rm \mathbf{u}}}^{h}\right\Vert _{L^{2}(\Omega)}$,
with ${\rm \mathbf{u}}$ the analytical solution given in Eq.\,\eqref{eq:exohdisp}
and ${\rm \mathbf{u}}^{h}$ the solution computed by the proposed
technique. The results show that all three types of elements have
monotonic convergence. For higher order elements, more accurate results
with similar number of DOF are obtained and the convergence rate is
also faster.

There are 37 out of 100 cells calculated for the stiffness matrices.
Among those 37 cells, 9 are regular quadtree cells and 28 are polygon
cells surrounding the hole. For the remaining cells, their stiffness
matrices are simply extracted from those 9 regular quadtree cells. 

\begin{figure}
\begin{centering}
\includegraphics[scale=0.8]{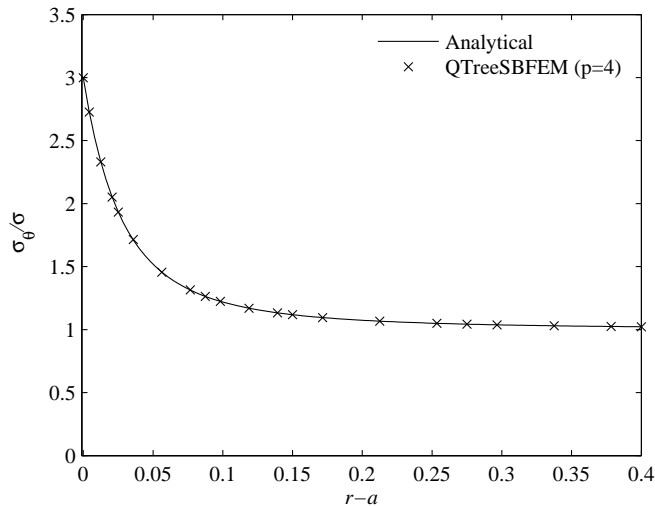}
\par\end{centering}

\caption{Thin square plate with a single circular hole under uniaxial tension}
\label{openholesigmatheta}
\end{figure}

To further demonstrate the accuracy of the proposed technique, $\sigma_{\theta}/\sigma$
along $A-B$ (see Fig.\,\ref{openhole}) is plotted in Fig.\,\ref{openholesigmatheta}
using the mesh given in Fig.\,\ref{ohm1} with 4th order elements.
It can be seen that the results of the proposed technique agree well
with the analytical solution, which has $\sigma_{\theta}/\sigma=3$
at $A$ ($\theta=90^{\circ},r=a$). For points away from $A$, $\sigma_{\theta}/\sigma$
approaches 1. 

\subsubsection{Approximation of infinite plate by varying $L/a$ ratio }

\begin{figure}
\begin{centering}
\includegraphics[scale=0.8]{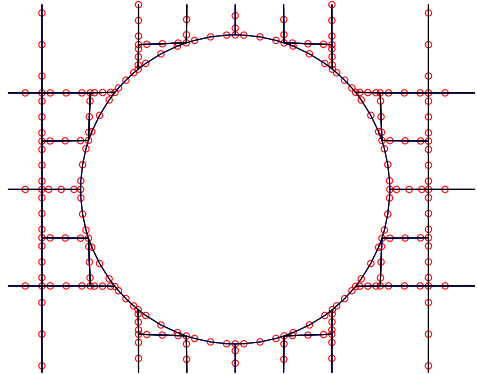}
\par\end{centering}

\caption{Cell pattern around the hole for $L/a=640$}
\label{ohm3}
\end{figure}

The same infinite plate can be approximated by increasing the $L/a$
ratio. The application of quadtree mesh facilitates such a study.
Only the left and right sides of the plate are subjected to uniaxial
in-plane tension stress $\sigma$. The element order used in this
study is $p=4$. The same mesh given in Fig.\,\ref{ohm1} is used
for $L/a=10$. The adaptive capability of quadtree mesh leads to the
same mesh pattern for all $L/a$ ratios. Fig.\,\ref{ohm3} shows
the cells around the hole for $L/a=640$ and it is exactly the same
as the one shown in Fig.\,\ref{ohm2}. 

For $L/a=640$, although there are 316 cells in total, only 37 cells
are calculated for the stiffness matrices, which is the same as the
previous study. The results of $\sigma_{\theta}/\sigma$ at $A$ with
varying $L/a$ ratio are given in Table\,\ref{openholesigmathetavsratio}.
It is seen that the analytical solution ($\sigma_{\theta}/\sigma=3$)
is quickly approached when increasing the $L/a$ ratio.

\begin{table}
\caption{Normalised stress ($\sigma_{\theta}/\sigma$) at $A$ of the thin
square plate with a circular hole}

\centering{}%
\begin{tabular}{ccc}
\hline 
$L/a$ ratio & No. of Nodes & $\sigma_{\theta}/\sigma$ at $A$\tabularnewline
\hline 
10 & 860 & 3.3591\tabularnewline
40 & 1428 & 3.0204\tabularnewline
160 & 1996 & 3.0049\tabularnewline
640 & 2564 & 2.9991\tabularnewline
\hline 
\end{tabular}\label{openholesigmathetavsratio}
\end{table}

\subsection{Square plate with multiple holes}

\begin{figure}
\begin{centering}
\includegraphics[scale=0.7]{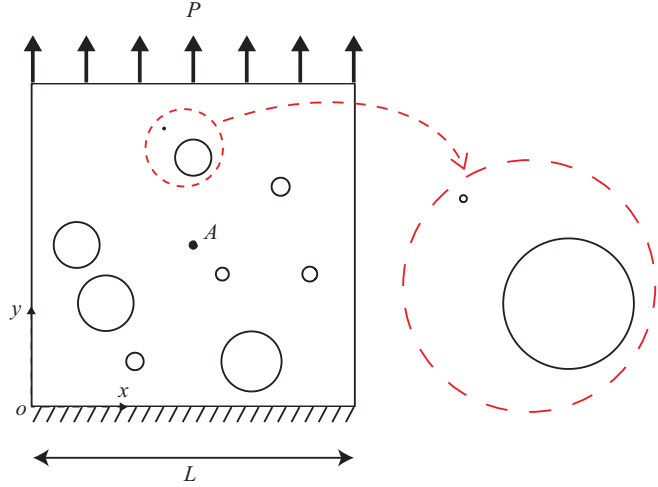}
\par\end{centering}

\caption{Thin square plate multiple holes under uniaxial tension}
\label{multiholes}
\end{figure}

A unit square plate with 9 randomly distributed holes of different
sizes, shown in Fig.\,\ref{multiholes}, is analysed. This example
highlights the automation and flexibility of the proposed technique
in handling the mesh transition between features with various dimensions.
The ability of capturing curved boundaries accurately using high-order
elements is also demonstrated. The displacements at the bottom edge
of the plate are fully constrained and a uniform tension $P=1$ is
applied at the top edge of the plate. A set of consistent units are
chosen.

To validate the results, the same problem is solved using the commercial
FEA software ANSYS V14.5. The plate is discretised using 8-node quadrilateral
elements (PLANE183). In order to demonstrate the automation and performance
of the proposed technique, similar user inputs are given to ANSYS
to generate a mesh for comparison. 

In ANSYS, the square plate is divided into 4 equal-sized quadrants
such that a centre key point is created for result comparison. The
9 holes are introduced to the plate by subtracting their areas from
the square plate. The mesh constructed in ANSYS is unstructured (paving).
A single variable $(N)$ is used to control the mesh density and is
used for mesh refinement. For each hole, the number of element around
them is equal to $4N$. And for all the straight lines, the size of
the elements is set to be $1/3N$, which gives each outer edge approximately
$3N$ elements. Note that the mesh used in ANSYS is far from optimal
and structured mesh should be used for better performance. However,
the main objective using paving mesh and only controlling the boundary
element divisions is to show how the proposed technique and ANSYS
perform when minimal number of controlling variables are used for
the meshing. For the proposed technique, seed points with $s_{b}=(4N\times s_{max})$,
where $s_{max}$ is the maximum allowed number of seed points in a
cell, are set on the circular holes to generate a mesh with similar
number of boundary division as the one in ANSYS. 

\begin{figure}
\begin{centering}
\subfigure[ANSYS (1358 elements)]{\label{mhm1}\includegraphics[scale=0.5]{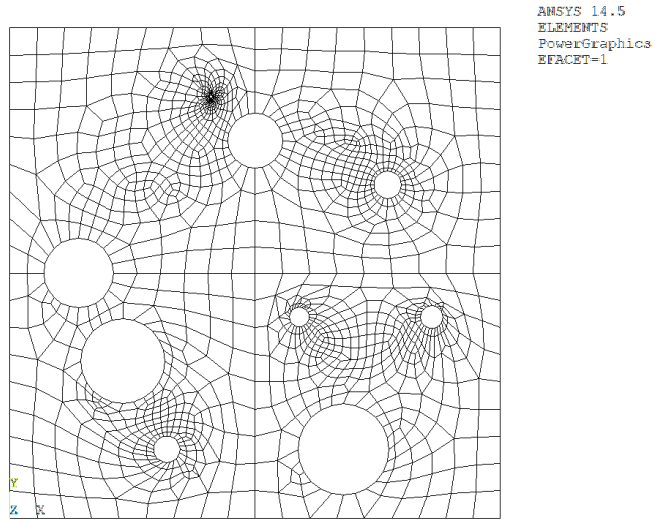}}
\par\end{centering}

\begin{centering}
\subfigure[QTreeSBFEM (1169 cells)]{\label{mhm2}\includegraphics{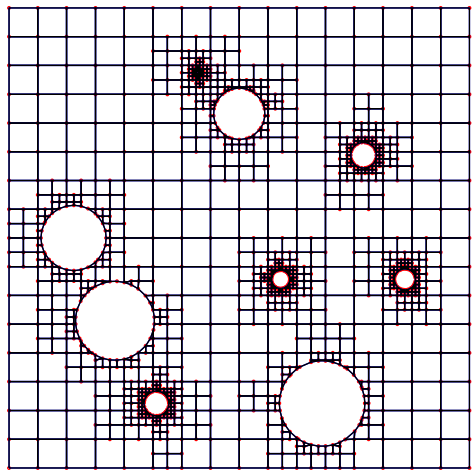}}
\par\end{centering}

\caption{Mesh of a square plate with 9 circular holes}
\label{multiholesmesh}
\end{figure}

Fig.\,\ref{multiholesmesh} shows the ANSYS mesh (1358 elements)
and the quadtree mesh (1169 cells). Both meshes can effectively handle
the mesh transition between the holes. As commented earlier, while
the mesh in ANSYS can be further improved by designing a structured
layout, it would also require additional time and investigation effort.
The amount of additional effort depends on the complexity of the geometry
and user experience. And for the proposed technique, the resulting
mesh is always in a structured manner (see Fig.\,\ref{mhm2}) without
additional effort. The time of generating the quadtree mesh is around
3s using the computer with details outlined in the beginning of this
section. 

For the stiffness calculation, there are 364 out of 1169 cells calculated.
Among those 364 cells, 12 are master quadtree cells and 352 are polygon
cells surrounding the holes. The stiffness matrices for all the other
cells are simply extracted from the calculated master quadtree cells.
The total time from constructing the stiffness to obtaining the displacement
solutions is less than 3.2s when using 5th order elements. 

\begin{figure}
\begin{centering}
\includegraphics[scale=0.8]{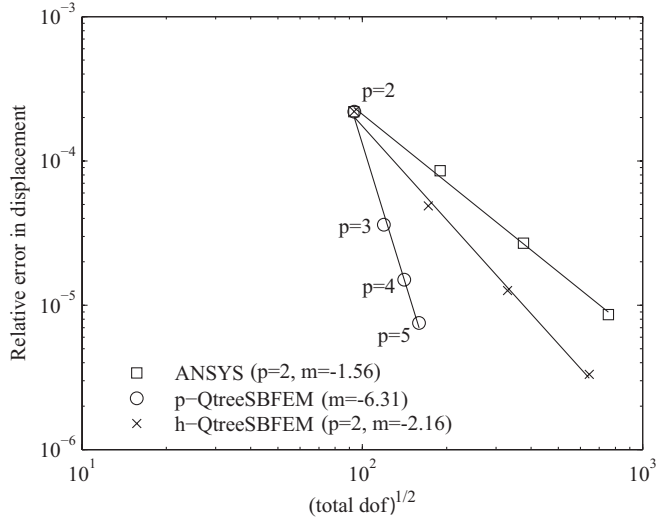}
\par\end{centering}

\caption{Convergence results of the centre displacement vector sum for the
square plate with multiple holes, where $p$ is the element order
and $m$ is the slope of the fitted line }
\label{multiholescon}
\end{figure}

Table\,\ref{multiholestab1} shows the convergence of the displacement
components at the centre point $A$ with increasing element order
using mesh in Fig.\,\ref{mhm2}. In order to highlight the convergence
performance, Fig.\,\ref{multiholescon} shows the relative error
of the present results of the displacement vector sum at point $A$.
The error is calculated based on the converged ANSYS results, which
converged to the first 6 significant digits. Also shown in the same
figure are the relative errors of the ANSYS results and another set
of results of the proposed technique. They are both generated through
a series of $h-$refinement with the use of 2nd order elements. It
is observed that the present results with $p-$refinement converge
with the fastest rate. And for the $h-$refinement, the present results
are basically converging at the same rate as those of ANSYS with slightly
better accuracy. The present results demonstrate that for the same
accuracy, much less number of DOFs is required when using high-order
elements, and curved boundaries are also accurately modelled with
minimal number of cells. 

\begin{table}
\caption{Centre displacement results of a square plate with 9 circular holes}

\centering{}%
\begin{tabular}{cccc}
\hline 
Elem. Order & No. of Nodes & $u_{x}\times10^{4}$ at $A$ & $u_{y}\times10^{3}$ at $A$\tabularnewline
\hline 
2 & 4379 & 4.98298 & 6.67236\tabularnewline
3 & 7157 & 4.98350 & 6.67359\tabularnewline
4 & 9935 & 4.98241 & 6.67374\tabularnewline
5 & 12713 & 4.98197 & 6.67379\tabularnewline
\hline 
\end{tabular}\label{multiholestab1}
\end{table}

\begin{figure}
\begin{centering}
\subfigure[ANSYS (285312 nodes)]{\label{mhsya}\includegraphics[scale=0.5]{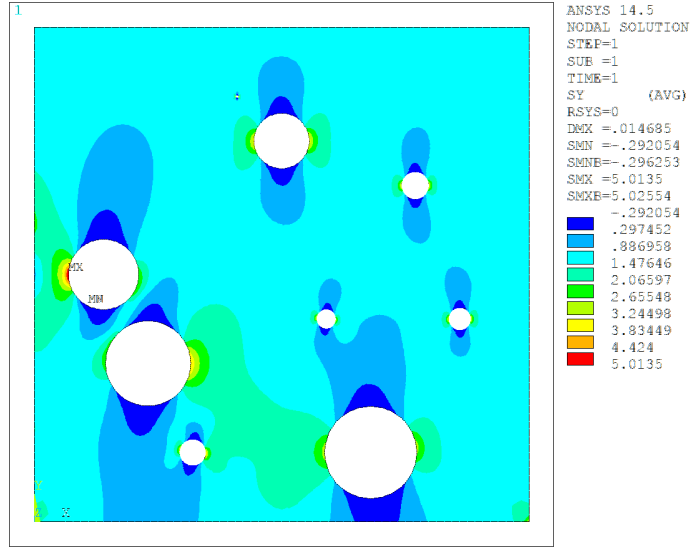}}
\par\end{centering}

\begin{centering}
\subfigure[QTreeSBFEM with 5th order elements (12713 nodes)]{\label{mhsy}\includegraphics{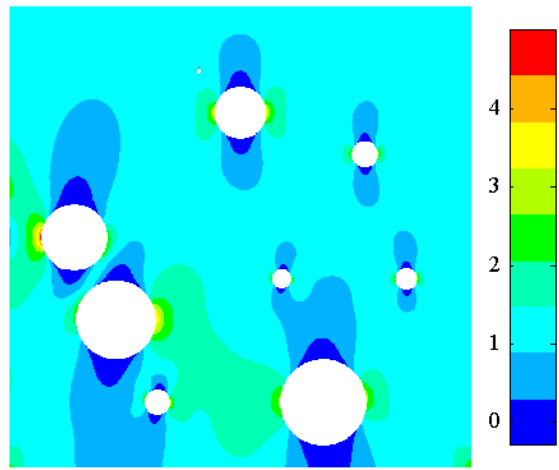}}
\par\end{centering}

\caption{$\sigma_{y}$ of the square plate with random pattern of holes}
\label{multiholessy}
\end{figure}

In order to further demonstrate the overall consistency of the present
results, Fig.\,\ref{multiholessy} shows the contour plots of $\sigma_{y}$,
from both ANSYS and the proposed technique. Good agreement is observed
from the contour plots.

\subsection{Square plate with multiple cracks emanating from a hole}

\begin{figure}
\begin{centering}
\includegraphics[scale=0.7]{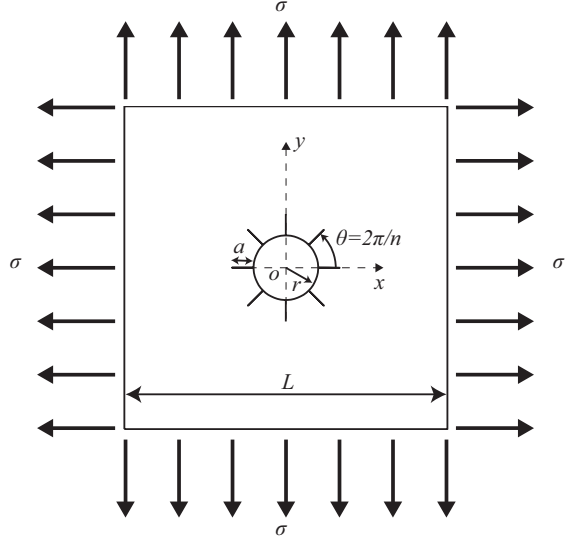}
\par\end{centering}

\caption{Thin square plate with cracks emanating from a hole under bi-axial
tension}
\label{holemulticracks}
\end{figure}

A square plate of length $L$ with a centre hole of radius of $r$
given in Fig.~\ref{holemulticracks} is considered. $n$ cracks with
crack length $a$ emanate from the hole. This example aims to show
the simplicity and effectiveness of the proposed technique to solve
problems with singularities.

In this example, to approach the assumption of an infinite plate,
$r/L=0.01$ is considered. A parametric study is performed considering
$n=2,4,8$ cracks surrounding the hole with various $s=\frac{a}{a+r}$
ratio. The element order used in this study is $p=4$, which is capable
to model the circular boundary accurately as shown in the first example.
The present results are compared with the reference solution of the
stress intensity factor given in \citet{tada2000stress}. 

\begin{figure}
\begin{centering}
\subfigure[4 cracks]{\label{hmcm1}\includegraphics{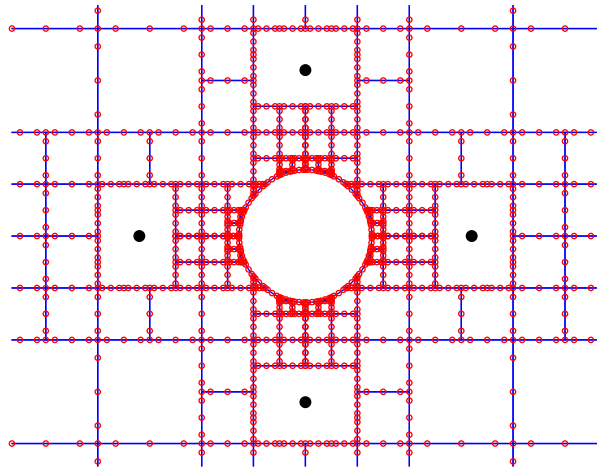}}
\par\end{centering}

\begin{centering}
\subfigure[8 cracks]{\label{hmcm2}\includegraphics{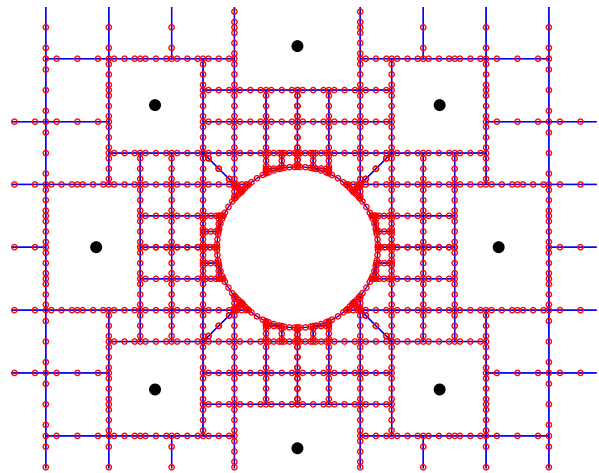}}
\par\end{centering}

\caption{Mesh of the square plate with cracks emanating from a hole. The black
dots represent the crack tips.}
\label{holemulticracksmesh}
\end{figure}

Fig.\,\ref{holemulticracksmesh} shows the mesh around the central
hole, with 4 and 8 cracks around the edge and $s=0.6$. Based on the
proposed technique, no refinement is required around the crack tips.
This facilitates the study with multiple cracks and results in less
computational effort when comparing to the conventional FEM. 

Fig.\,\pageref{holemulticracksres} shows the stress intensity factor
($F_{I}=K_{I}/(\sigma\sqrt{\pi a})$) computed from the proposed technique
for different value of $s$. Excellent agreement with the reference
solution \citep{tada2000stress} is observed. This demonstrates the
accuracy of the proposed technique in dealing with stress singularities
as well as the feasibility in handling geometry with complicated features. 

\begin{figure}
\begin{centering}
\includegraphics[scale=0.8]{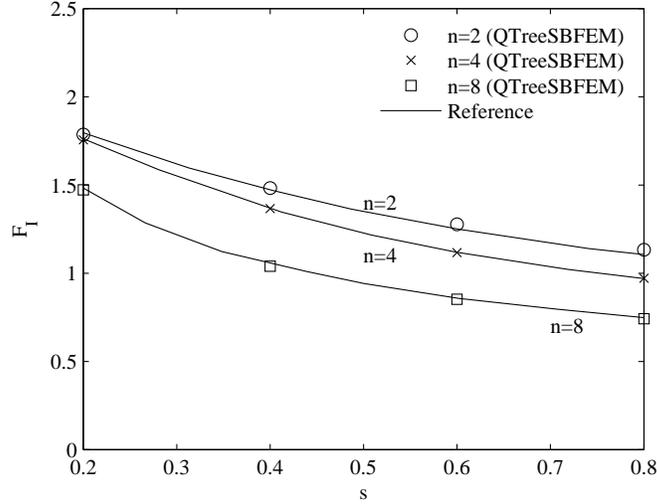}
\par\end{centering}

\caption{Stress intensity of the square plate with cracks emanating from a
hole}
\label{holemulticracksres}
\end{figure}

\subsection{Square plate with two cracks cross each other }

\begin{figure}
\begin{centering}
\includegraphics[scale=0.7]{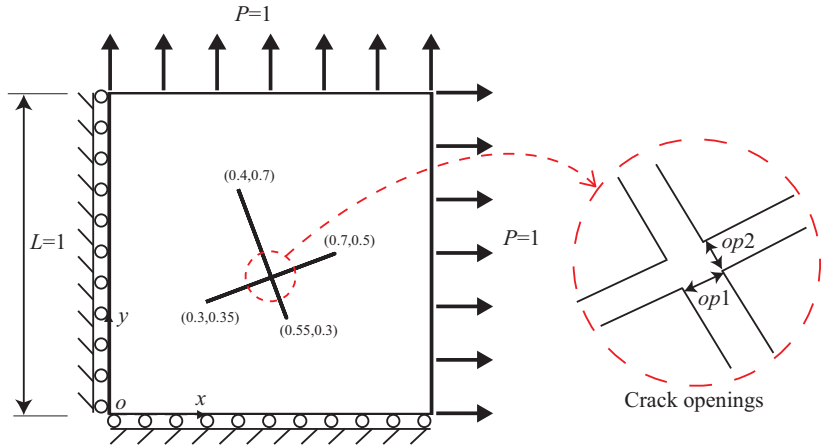}
\par\end{centering}

\caption{Thin square plate with two cracks cross each other}
\label{crosscrack}
\end{figure}

A square plate of length $L$ with two cracks cross each other is
considered. The dimensions of the plate and the cracks as well as
the boundary conditions are shown in Fig.~\ref{crosscrack}. This
example highlights the automatic mesh generation of the proposed technique
and the capacity to handle problems with complicated crack configuration. 

\begin{figure}
\begin{centering}
\includegraphics{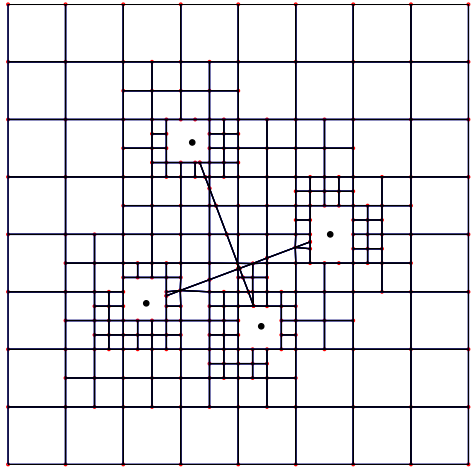}
\par\end{centering}

\caption{Mesh of a square plate with two cracks cross each other (216 cells)}
\label{crosscrackmesh}
\end{figure}

Fig.\,\ref{crosscrackmesh} shows the quadtree mesh of the proposed
technique. The mesh only requires defining seed points on the domain
boundary, along the cracks and around the crack tips to control the
quadtree mesh density. The mesh generation is fully automatic without
the requirement of dividing area regions. The resulting mesh contains
a total of 216 cells. During the construction of the stiffness matrix,
only 48 cells are computed, which contains 9 master quatree cells
and 39 polygon cells. 

For the same problem, it would require a few more steps to generate
a mesh in FEA. These include defining crack tip regions that directly
affect the solution accuracy, and designing proper refinement strategy
that directly affects the convergence performance. For example in
ANSYS, a command \emph{``}KSCON'' needs to be issued to each crack
tip in order to generate two circular layers of elements (1 layer
singular elements) around the tip. The radius of the two circular
layers of elements is solely based on user experience and \emph{trial-and-error}.
Shape warning on the elements would occur if the settings of that
command are not consistent with the global mesh. Moreover, automatic
$h-$refinement is not applicable when ``KSCON'' is activated. It
would, therefore, require multiple steps to conduct convergence study
with $h-$refinement, such as reducing the radius of the circular
layers of element around the crack tips and increasing elements in
circumferential direction around the crack tips.

Table\,\ref{crosscracktab} shows the crack opening displacements
($op1$ and $op2$ in Fig.\,\ref{crosscrack}) with increasing element
order. Similar to previous examples, the present results converge
rapidly with minimal number of nodes increased. The results between
using the 2nd order elements and using the 5th order elements are
different with less than $0.04\%$. 

\begin{table*}
\caption{Crack opening displacement: $u_{op1}$ for the opening $op1$ and
$u_{op2}$ for the opening $op2$}

\noindent \centering{}%
\begin{tabular}{cccc}
\hline 
Elem. Order & No. of Nodes & $u_{op1}\times10^{-3}$ & \multicolumn{1}{c}{$u_{op2}\times10^{-3}$}\tabularnewline
\hline 
2 & 818 & 5.1300 & 6.9710\tabularnewline
3 & 1333 & 5.1274 & 6.9727\tabularnewline
4 & 1848 & 5.1279 & 6.9726\tabularnewline
5 & 2363 & 5.1280 & 6.9725\tabularnewline
\hline 
\end{tabular}\label{crosscracktab}
\end{table*}

\subsection{Cracked nuclear reactor under internal pressure}

\begin{figure}
\begin{centering}
\includegraphics{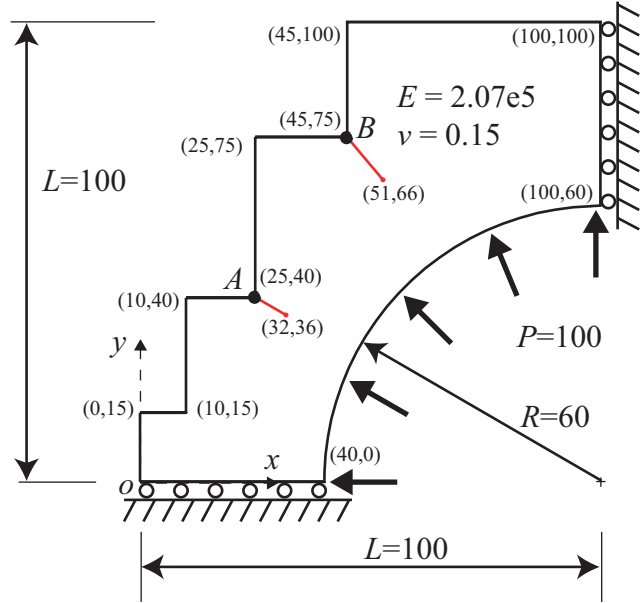}
\par\end{centering}

\caption{Cracked nuclear reactor under internal pressure}
\label{nreact}
\end{figure}

In this final example, a nuclear reactor under internal pressure \citep{Simpson2013}
is analysed. Due to symmetry, only a quadrant of the reactor is modelled.
The geometry, material properties, loading and dimension are shown
in Fig.\,\ref{nreact}. Also shown in the figure are the two cracks
introduced on the outer boundary. This example shows the flexibility
of the proposed technique and the developed meshing algorithm to model
more practical structures. 

\begin{figure}
\begin{centering}
\includegraphics{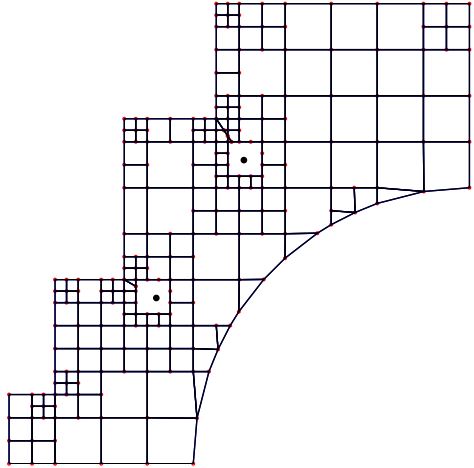}
\par\end{centering}

\caption{Mesh for the quadrant of the reactor with 2 cracks (160 cells)}
\label{nreactmesh}
\end{figure}

Fig.\,\ref{nreactmesh} shows the quadtree mesh used in this example,
which contains a total of 160 cells. Using the proposed technique
only requires seed points to be defined at the boundaries to control
the quadtree mesh density. No additional requirement for the cells
containing the crack tips is necessary. The time spent on generating
the quadtree mesh is less than 0.8s using the same computer with details
outlined in the beginning of this section. 

The calculation of stiffness matrix involves computing 44 out of the
160 cells, in which 12 are master quadtree cells and 32 are polygon
cells. The total time from constructing the stiffness matrix to obtaining
the displacement solution is less than 0.7s when using 5th order elements. 

A convergence study is conducted by increasing the order of the element
without changing the quadtree layout in Fig.\,\ref{nreactmesh}.
Table\,\ref{nreacttab1} shows the two crack opening displacements
at points $A$ and $B$ (Fig.\,\ref{nreact}). The present results
converge quickly with the element order increased. The difference
between using the 2nd order elements and using the 5th order elements
is less than 0.1\%. This further highlights the advantage of using
high-order elements that can model curved boundary more accurately
with minimal number of cells. 

\begin{table*}
\caption{Crack opening displacement: $u_{A}$ for opening at$(25,40)$, $u_{B}$
for opening at $(45,75)$}

\noindent \centering{}%
\begin{tabular}{cccc}
\hline 
Elem. Order & No. of Nodes & $u_{A}\times10^{2}$ & \multicolumn{1}{c}{$u_{B}\times10^{2}$}\tabularnewline
\hline 
2 & 635 & 7.14372 & 2.60404\tabularnewline
3 & 1031 & 7.15066 & 2.60336\tabularnewline
4 & 1427 & 7.15077 & 2.60334\tabularnewline
5 & 1823 & 7.15077 & 2.60332\tabularnewline
\hline 
\end{tabular}\label{nreacttab1}
\end{table*}

\section{Conclusion\label{sec:Conclusion}}

This paper has presented a numerical technique to automate stress
and fracture analysis using the SBFEM and quadtree mesh of high-order
elements. Owing to the nature of the SBFEM, the proposed technique
has no specific requirement, such as deriving conforming shape functions
or sub-triangulation, to handle quadtree cells with hanging nodes.
High-order elements are used within each quadtree cell directly. 

The quadtree mesh generation is fully automatic and involves minimal
number of user inputs and operation steps. Boundaries are modelled
with scaled boundary polygons and this allows the proposed technique
to conform the boundary without excessive mesh refinement. The meshing
algorithm is also applicable for problems with singularities. The
use of quadtree mesh leads to an efficient approach to compute the
global stiffness matrix. This facilitates the analysis that requires
a significant number of cells using high-order elements. Five numerical
examples are presented to highlight the functionality and performance
of the proposed technique. The present results show excellent agreement
with analytical solutions and those computed by the FEM.

\section*{Reference}

\section*{\textmd{\normalsize{\bibliographystyle{elsart-num-names}
\bibliography{Quadtree}
}}}
\end{document}